\documentclass[10pt]{article}
\usepackage{amssymb,amsmath}
\usepackage{epsf}
\usepackage{times,bm,mathrsfs}  % Fonts

\newcommand{\qedsymb}{\hfill{\rule{2mm}{2mm}}}
\newenvironment{proof}[1][]{\begin{trivlist}
\item[\hspace{\labelsep}{\bf\noindent Proof#1:\/}] }{\qedsymb\end{trivlist}}

\newcounter{AbcT}

\newtheorem {Theorem}    {Theorem}[section]

\newtheorem {Lemma}      [Theorem]    {Lemma}
\newtheorem {Corollary}  [Theorem]    {Corollary}

\newcommand{\ignore}[1]{}

\def\E{{\bf{E}}}
\def\P{{\bf{P}}}

\def\C{{\mathbb{C}}}

\def\v0{{\bf 0}}

\def\0{\hat{0}}
\def\1{\hat{1}}

\def\path{{\tt path}}

\def\lam{\lambda}

\def\phi{\varphi}

\def\be{\begin{equation}}
\def\ee{\end{equation}}

\def\TVD#1{\Vert{#1}\Vert_{\rm TV}}
\def\tmix{\tau_{\rm mix}}
\def\Unif{{\cal U}}

\def\wdk{\widetilde{D}_k}

\begin{document}
\title{Shuffling by semi-random transpositions}
%\small Draft \#2.11, October 14}
\author{Elchanan Mossel \thanks{Supported by
a Miller Fellowship in CS and Statistics, U.C. Berkeley.}
\\ Statistics \\ U.C. Berkeley \\ mossel@stat.berkeley.edu \and
Yuval Peres \thanks{Research supported in part by NSF Grants
DMS-0104073 and DMS-0244479. Part of this work was done while the
author was visiting Microsoft Research. } \\ Statistics and Mathematics\\
U.C. Berkeley \\ peres@stat.berkeley.edu \and
Alistair Sinclair \thanks{Supported in part by NSF grant CCR-0121555.
Part of this work was done while the author was on sabbatical leave at
Microsoft Research and at Ecole Polytechnique, Paris.}
 \\ Computer Science \\ U.C. Berkeley \\
sinclair@cs.berkeley.edu}

\maketitle

%%%%%%%%%%%%%%%%%%%%%%%%%%%%%%%%%%%%%%%%%%%%%%%%%%%%%%%%%%%%%%%%%%%%%%%
\begin{abstract}
%%%%%%%%%%%%%%%%%%%%%%%%%%%%%%%%%%%%%%%%%%%%%%%%%%%%%%%%%%%%%%%%%%%%%%%
  In the cyclic-to-random shuffle, we are given $n$ cards arranged in
  a circle. At step $k$, we exchange the $k$'th card along the circle
  with a uniformly chosen random card. The problem of determining the
  mixing time of the cyclic-to-random shuffle was raised by Aldous and
  Diaconis in 1986. Recently, Mironov used this shuffle as a model for
  the cryptographic system known as ``RC4'' and proved an upper bound
  of $O(n \log n)$ for the mixing time.  We prove a matching lower bound,
  thus establishing that the mixing time is indeed of order $\Theta(n \log n)$.
  We also prove an upper bound of $O(n\log n)$ for the mixing time of
  any ``semi-random transposition shuffle'', i.e., any shuffle in which
  a random card is exchanged with another card chosen according to an
  arbitrary (deterministic or random) rule.  To prove our lower bound,
  we exhibit an explicit
  complex-valued test function which typically takes very different
  values for permutations arising from the cyclic-to-random-shuffle
  and for uniform random permutations; we expect that this test function may be
  useful in future analysis of RC4.  Perhaps surprisingly, the
  proof hinges on the fact that the function $e^z-1$ has nonzero fixed
  points in the complex plane. A key insight from our work is the
  importance of complex analysis tools for uncovering structure in
  nonreversible Markov chains.

%%%%%%%%%%%%%%%%%%%%%%%%%%%%%%%%%%%%%%%%%%%%%%%%%%%%%%%%%%%%%%%%%%%%%%%%%%%
\end{abstract}
%%%%%%%%%%%%%%%%%%%%%%%%%%%%%%%%%%%%%%%%%%%%%%%%%%%%%%%%%%%%%%%%%%%%%%%%%%%

%%%%%%%%%%%%%%%%%%%%%%%%%%%%%%%%%%%%%%%%%%%%%%%%%%%%%%%%%%%%%%%%%%%%%%%%%%%
%%%%%%%%%%%%%%%%%%%%%%%%%%%%%%%%%%%%%%%%%%%%%%%%%%%%%%%%%%%%%%%%%%%%%%%%%%%
\section{Introduction}
\label{sec:intro}
%%%%%%%%%%%%%%%%%%%%%%%%%%%%%%%%%%%%%%%%%%%%%%%%%%%%%%%%%%%%%%%%%%%%%%%%%%%
%%%%%%%%%%%%%%%%%%%%%%%%%%%%%%%%%%%%%%%%%%%%%%%%%%%%%%%%%%%%%%%%%%%%%%%%%%%

The {\it mixing time\/} of a Markov chain on a finite state space is the
number of steps until it is close to its stationary
distribution, starting from an arbitrary state.  The mixing time is a
key parameter in analyzing random sampling algorithms
and is of intrinsic interest in probability and statistical physics as well.
For many natural Markov chains,
if some of the randomness is removed from the transition rule,
resulting in a ``more deterministic'' process with the same stationary
distribution, the chain becomes significantly harder to analyze.
Indeed,  some of the most
challenging problems in the field concern the analysis of such
``pseudo-random'' variants of well understood chains.
Some examples include the riffle shuffle~\cite{Gilbert,Reeds}
compared to the thorp shuffle~\cite{Thorp},
the asymmetric exclusion process~\cite{BBHM} compared with
its systematic scan version~\cite{DR} and the comparison between the
standard and systematic versions
of Glauber dynamics for Gaussian fields \cite{GS,Amit1,Amit2,AG}.

%{\bf Cyclic-to-random shuffle.}
Shuffling by random transpositions is one of the simplest random walks on the
symmetric group: given $n$ cards
in a row, at each step two cards are picked uniformly at random and
exchanged. This shuffle was precisely analyzed in 1981, see~\cite{DiSh}.
In the ``cyclic-to-random'' shuffle
(invented by Thorp \cite{Thorp2}), at step~$t$
a uniformly chosen  random card is exchanged with the card at position
$\, t \! \mod n$.
It is easy to see that this semi-random shuffle still converges
to the uniform distribution on permutations of $n$ cards.
In their landmark 1986 paper on card shuffling~\cite{AlDi},
Aldous and Diaconis
posed as a challenge the analysis of the cyclic-to-random shuffle.  More recently,
Mironov~\cite{Mironov} related this shuffle to the behavior of
the RC4 encryption algorithm. Mironov showed that a strong uniform time argument,
due to Broder, can be adapted to yield an upper bound of $O(n\log n)$
on the mixing time.  He posed as an open problem 
whether this bound is tight.

In this paper we establish a lower
bound of $\Omega(n\log n)$ for the mixing time of the cyclic-to-random shuffle,
 thus answering the questions posed by
Aldous and Diaconis and by Mironov.    We also prove a general upper
bound of $O(n\log n)$ on the mixing time of {\it any\/} semi-random
 transposition shuffle, i.e., any shuffle in which
a random card is exchanged with another card chosen according
to an arbitrary (deterministic or random) rule that may vary
at each step.  Previously, the best available upper bound for such a general
process was $O(n^2)$,  proved by Pak~\cite{Pak}.

To prove the lower bound for the cyclic-to-random transposition shuffle
$\{\sigma_t\}$, we  find an eigenfunction of the shuffle that mixes slowly.
(This approach was used by Wilson~\cite{wilson1,wilson2} to prove
$\Omega(n^3 \log n)$ lower bounds for the shuffle generated by transpositions
of  adjacent cards and several variants.)  First, we
determine the eigenvalues of a nonreversible renewal
Markov chain $M$ on the
$n$-cycle which describes the behavior of a single card. The asymptotics for the
leading eigenvalues of $M$ depend on the fact that the function
$e^z-1$ has nonzero fixed points in the complex plane.
We then pick an eigenfunction $f$ for $M$ and use it to construct a
test function $F$, defined on permutations, which is a  weighted sum of $f$
applied to the locations of all cards.
To show that the distribution at time $t$ of $F(\sigma_t)$
is far from the distribution of $F(\sigma)$ for a uniform  random
permutation $\sigma$, the key is to estimate the variance.
The variance is a sum of correlations between pairs of cards;
to bound these correlations, we couple the shuffle with a system of independent
particles evolving according to $M$. This coupling approach
has intuitive appeal, and could potentially be used for other chains on permutations.
Alternatively, one could bound the variance of $F(\sigma_t)$
using the martingale decomposition method
of Wilson~\cite{wilson1,wilson2}.

Our general upper bound for semi-random transpositions is proved via a strong uniform
time argument, extending earlier arguments of Broder and Mironov.

We believe that some of our technical insights may be carried over to other
situations where lower bounds for nonreversible or
``pseudo-random'' Markov chains are
sought. These insights include:
\begin{itemize}
\item
The analysis of a given Markov chain with a transition rule that
varies in time can sometimes be reduced to the analysis of an
equivalent time-homogeneous chain.
%This is crucial in simplifying the analysis.
\item
Coupling arguments, which are often applied to obtain upper bounds
 for mixing times, can also be used to establish lower
bounds.
\item
When seeking to understand a nonreversible Markov chain,
 results of classical complex analysis such as
Rouch\'e's theorem, can be powerful tools.
Thus methods from complex analysis should be added to
 techniques from probability,
combinatorics, functional analysis and representation theory
in the toolkit of Markov chain analysis.

\end{itemize}

\subsection{Statement of main results}

  Let
  $\{L_t\}_{t=1}^\infty$ be a  sequence of random variables taking values
   in $[n]=\{0,1,\dots,n-1\}$ and let $\{R_t\}_{t=1}^\infty$
    be a sequence of i.i.d.\ cards chosen uniformly from $[n]$.
       The {\bf semi-random
  transposition shuffle\/} generated by~$\{L_t\}$ is a stochastic
  process $\{\sigma_t^*\}_{t=0}^\infty$ on  the symmetric group
$S_n$, defined as
  follows.  Fix the initial permutation~$\sigma^*_0$.
  The permutation~$\sigma_t^*$ at time~$t$ is obtained from~$\sigma^*_{t-1}$
  by transposing the cards at locations $L_t$ and $R_t$.

In the {\bf cyclic-to-random shuffle}, the sequence $L_t$ is given by
$L_t = t \! \mod n$.

The stochastic process~$\{\sigma_t^*\}$ is a time-inhomogeneous
Markov chain on~$S_n$, and converges to the uniform stationary
distribution for any~$\sigma_0^*$ and any choice of~$\{L_t\}$.
It is a time-homogeneous Markov chain if the~$L_t$ are i.i.d.\
The special case where the $L_t$ are i.i.d.\ uniform  is
the  random transposition
shuffle~\cite{Aldous,AlDi, DiSh}, the random walk on~$S_n$ generated
by all transpositions; at the other extreme, if all the~$L_t$
are identically 0, we get the random
walk generated by ``star transpositions'', where in each step  a randomly
chosen card is exchanged with the card in position 0.
% This is also called the ``top-with-random transposition shuffle.''

Let  $\mu_t^*$ be the distribution of $\sigma_t^*$
at time~$t$, and let $\TVD{\mu_t^*-\Unif}$
denote the total variation distance between  $\mu_t^*$
and the uniform distribution~$\Unif$. Define the {\it mixing time\/} by $$
   \tmix = \max_{\sigma_0}\min\{t:\TVD{\mu_t^*-\Unif}\le\frac{1}{2 e}\}.$$
The choice of constant~$\frac{1}{2e}$ ensures that, for any $\epsilon>0$,
we have $\TVD{\mu_t^*-\Unif}\le\epsilon$ if
$t\ge \lceil\log\epsilon^{-1}\rceil\tmix$ (see~\cite{Aldous}).

\begin{Theorem} \label{th:lower}
The cyclic-to-random transposition shuffle has mixing time $\Omega(n\log n)$.
More precisely, the mixing time is at least
\[
\frac{n \log n}{|\zeta+1|(1 + o(1))},
\]
where $\zeta$ is any nonzero complex root of the  equation
$\psi(z) =e^z - z - 1 = 0$.
\end{Theorem}
Using Mathematica, we find the root $\zeta=
2.088... + 7.461... \times i$ of $\psi$. This gives
%2.0888430156... + 7.4614892856... \times i$ of $\psi$. This gives
$|1+\zeta| = 8.075...$ and yields a lower bound of
$(.123 + o(1)) n \log n$ for the mixing time.
%$|1+\zeta| = 8.0755664528...$ and yields a lower bound of
%$(.1238303227... + o(1)) n \log n$ for the mixing time.

\begin{Theorem}
\label{th:upper}
The semi-random transposition shuffle $\{\sigma_t^*\}$
generated by any sequence~$\{L_t\}$,
has mixing time at most $O(n\log n)$.
More precisely, there is a constant $C_0$
such that for any $C_1 > C_0$ and any
 initial configuration $\sigma_0^*$, we have
\[
\TVD{\mu_t^* - \Unif} \leq n^{-\beta} \mbox{ \rm for all }
t>C_1 n \log n \,,
\]
for some $\beta=\beta(C_1)>0$.
%(Recall that $\mu_t^*$ is the distribution of $\sigma_t^*$
%at time $t$.)
\end{Theorem}
{\bf Remark.} The proof shows that we
 can take $C_0= 32\theta^{-3}+\theta^{-1}$ where
$\theta =  e^{-2}(1-e^{-1})/2$. We do not know the minimal value of $C_0$;
it cannot be strictly less than 1 because of
the star transpositions
shuffle, where the mixing time is $(1+o(1)) n \log n$,
see \cite{diaconis2}.
%%%%%%%%%%%%%%%%%%%%%%%%%%%%%%%%%%%%%%%%%%%%%%%%%%%%%%%%%%%%%%%%%%%%%%%%%%%
\section{A lower bound for the cyclic-to-random shuffle}
%%%%%%%%%%%%%%%%%%%%%%%%%%%%%%%%%%%%%%%%%%%%%%%%%%%%%%%%%%%%%%%%%%%%%%%%%%%
%%%%%%%%%%%%%%%%%%%%%%%%%%%%%%%%%%%%%%%%%%%%%%%%%%%%%%%%%%%%%%%%%%%%%%%%%%%
\subsection{The behavior of a single card via renewals}
Fixing a specific card $a$,
it is natural to study the renewal chain on the
state space $[n]=\{0,\ldots,n-1\}$,
where state $i \in [n]$ indicates that the location
$j$ of card $a$ satisfies $j + i = t \mod n$.  This chain is
described by the transition matrix~$M$, where for all $i \in [n]$, we have
$M_{0,i} = 1/n$ and $M_{i,1} = 1/n$, while $M_{i, i+1} = 1 - 1/n$,
for all $i \geq 1$. (For $i=n-1$, the last equation
reads $M_{n-1,0} = 1 - 1/n$.)
In other words,
\begin{equation*}  %\label{eq:defM}
M = \left( \begin{array}{ccccccc}
\frac{1}{n} & \frac{1}{n} & \frac{1}{n} & \frac{1}{n} & \ldots
& \frac{1}{n} & \frac{1}{n}\\
0 & \frac{1}{n} & 1 - \frac{1}{n} & 0 & 0 & \ldots & 0 \\
0 & \vdots & 0 & \ddots & 0 & \ldots & 0 \\
0 & \vdots & 0 & 0 & \ddots & \ldots & 0 \\
0 & \vdots & 0 & \ldots & 0 & \ddots & 0 \\
0 & \frac{1}{n} & 0 & \ldots & 0 & 0 & 1 - \frac{1}{n} \\
1 - \frac{1}{n} & \frac{1}{n} & 0 & \ldots & 0 & 0 & 0
\end{array} \right).
\end{equation*}

We will now find the eigenfunctions of the chain, that is, the right
eigenvectors of the matrix~$M$. Let $f = (f_0,\ldots,f_{n-1})^T$
 be such a (column) eigenvector.
Then we obtain the following equations:
\begin{equation}  \label{esum}
\frac{1}{n} \sum_{j=0}^{n-1} f_j = \lam f_0,
\end{equation}
and, for $1 \leq i \leq n-1$,
\begin{equation} \label{eq:eigeni}
\frac{1}{n} f_1 + (1 - \frac{1}{n}) f_{i+1} = \lam f_i \,\,\,\,\,
%(\mbox{ indices are computed modulo } n).
(\mbox{\rm where we denote  } f_n=f_0).
\end{equation}
It is easy to check that, up to scaling, $(1,\ldots,1)^T$ is the unique
eigenvector corresponding to the eigenvalue $\lam = 1$, and
that $(-1,n-1,-1,\ldots,-1)^T$
is the unique eigenvector corresponding to $\lam = 0$.

We now assume that $f$ is a right eigenvector corresponding
to an eigenvalue $\lam \notin \{0,1\}$.  Since $M$
is doubly stochastic, (\ref{esum}) implies that $\sum_{i=0}^{n-1} f_i = 0$
and $f_0=0$; to verify this, sum
(\ref{esum}) and the $n-1$ equations in (\ref{eq:eigeni}).

Writing $y_i = f_{i+1} - f_i$ for $1 \leq i \leq n-1$ (recall that
$f_n=f_0$) the equation (\ref{eq:eigeni}) for $i=1$ gives:
\begin{equation*} %\label{eq:eigen1}
 y_1 = {\frac{n(\lam-1)}{n-1}}f_1
%\frac{1}{n} x_1 + (1 - \frac{1}{n}) x_2 = \lam x_1 \,.
\end{equation*}
For $1 \leq i \leq n-2$, subtracting successive equations in
(\ref{eq:eigeni}) yields
\begin{equation*} %\label{eq:eigen1'}
(1 - \frac{1}{n}) y_{i+1} = \lam y_i.
\end{equation*}
Thus if we set $\gamma = \frac{n\lam}{n-1}$, then
$y_1 = (\gamma-{\frac{n}{n-1}})f_1$ and
$y_j = \gamma^{j-1} y_1$ for $2 \le j \le n-1$.
Without loss of generality we may assume that $f_1 = 1$.
Therefore, % with the notation $f_n=f_0$, we have
\begin{equation} \label{eq:x_i}
 f_{k} = 1 + \sum_{j=1}^{k-1} y_j = 1 + y_1
\sum_{j=1}^{k-1} \gamma^{j-1} =1 + \left( \gamma - \frac{n}{n-1} \right)
                \sum_{j=1}^{k-1} \gamma^{j-1}
\end{equation}
for $1 \leq k \leq n$. Thus
\begin{equation} \label{eq:xsum}
 (n-1) (1-\gamma) f_{k} = \Bigl(n-(n-1)\gamma \Bigr)\gamma^{k-1}-1
\end{equation}
for $1 \leq k \leq n$.
Since $\sum_{k=0}^{n-1} f_i=0$ and $f_n=f_0$, we infer that
$$
0=(n-1)(1-\gamma)^2 \sum_{k=1}^n f_{k}=\Bigl(n-(n-1)\gamma \Bigr)(1-\gamma^{n})
-n(1-\gamma) = \gamma-n\gamma^n+(n-1)\gamma^{n+1} \,.
$$
Since $\gamma \ne 0$, it follows that
\begin{equation} \label{poly}
(n-1) \gamma^n - n \gamma^{n-1} + 1 \,.
\end{equation}
 Note that this equation has a double root at $\gamma = 1$.
  We therefore conclude that the eigenvalues
$\lam \neq 0,1$ correspond (via the relation $\gamma = \frac{n\lam}{n-1}$)
 to the roots $\gamma \neq 1$
of (\ref{poly}).  We investigate
these roots next.

\subsection{Properties of roots of equation (\ref{poly}).}
\label{sec:roots}

\begin{Lemma} \label{lem:delta_leq_1}
All the solutions of (\ref{poly}) satisfy $|\gamma| \leq 1$.
\end{Lemma}

\begin{proof}
If $|\gamma| > 1$, then
\[
|(n-1) \gamma^{n-1}| > \Bigl|\sum_{i=0}^{n-2} \gamma^i\Bigr| =
\Bigl|\frac{\gamma^{n-1}-1}{\gamma-1}\Bigr|,
\]
and multiplying by $|\gamma-1|$ gives
\[
|(n-1) \gamma^n - (n-1) \gamma^{n-1}| > |\gamma^{n-1}-1|,
\]
so $\gamma$ cannot be a solution to (\ref{poly}).
\end{proof}

In the other direction, we need to show that (\ref{poly})
has solutions close to $1$. We prove:
\begin{Lemma} \label{lem:delta_geq_1}
There exists a solution of the equation
$(n-1) \gamma^n - n \gamma^{n-1} + 1 = 0$ which satisfies
\begin{equation} \label{eq:gamma_close}
|1 - \gamma| \leq  \frac{|\zeta|}{n}+O(\frac{1}{n^2}),
\end{equation}
and
\begin{equation} \label{eq:lam_close}
|1 - \lam| \leq \frac{|\zeta+1|}{n}+O(\frac{1}{n^2}),
\end{equation}
where $\zeta$ is any nonzero root of $e^\zeta - \zeta - 1 = 0$,
and $\lam = (1 - 1/n) \gamma$.
\end{Lemma}

\begin{proof}
By defining $\omega = \gamma^{-1}$, we obtain the
equation $\omega^n - n \omega + n-1 = 0$, or
$\omega^n + n (1 - \omega) - 1 = 0$.
Now write $\omega = 1 + z/n$ to get the asymptotic equation
$\psi(z)=e^z - z - 1 = 0$.
By Hurwitz's theorem (see \cite{Ahlfors}), every solution $\zeta$ of the
equation $\psi(\zeta)=0$ is a limit of solutions $z_n$ of the
equations $(1 + z_n/n)^n - z_n -1= 0$.
Since $\omega - 1 = z_n/n$, we obtain
\begin{equation} \label{eq:gamma_close2}
\gamma = 1 - z_n/n + O(\frac{1}{n^2}).
\end{equation}
Therefore
\begin{equation} \label{eq:lam_close2}
\lam = (1-1/n)\gamma=1 - \frac{1+z_n}{n} + O(\frac{1}{n^2}) \,.
\end{equation}
To get more precise estimates, recall $\psi(z)=e^z - z - 1$
and let $\varphi_n(z)=(1 + z/n)^n-z-1$.
By Taylor expansion,
$$
|n\log(1 + z/n) -z| = \frac{|z|^2}{2n} +O(\frac{1}{n^2}) \,,
$$
so in a bounded domain,
$$
|\varphi_n(z)-\psi(z)|=|(1 + z/n)^n -e^z| = \frac{|z^2 e^z |}{2n}
+O(\frac{1}{n^2}) \,.$$

Below we will prove that the equation $e^z-z-1=0$ has nonzero
roots. Let $\zeta$ be such a root, then $\zeta$ is a simple root,
 since $\psi'(\zeta)=e^{\zeta}-1=\zeta$.
Thus for $z$ on the circle $\{|z-\zeta|=b/n\}$, we have
$$
|\psi(z)| = |\psi'(\zeta)| \frac{b}{n} +O(\frac{1}{n^2})=
|\zeta| \frac{b}{n} +O(\frac{1}{n^2}).
$$
On the other hand, for $z$ on that circle,
$$
|\varphi_n(z)-\psi(z)|= \frac{|\zeta^2 e^\zeta |}{2n} +O(n^{-2}) \,.
$$
By Rouch\'e's Theorem (see \cite{Ahlfors}), it follows that if
$b>|\zeta e^\zeta/2|$ and $n$ is
large enough, then $\varphi_n$ has the same number of zeros as $\psi$
in the disk $\{|z-\zeta|< b/n\}$, namely, exactly one zero.
We thus obtain (\ref{eq:gamma_close}) by (\ref{eq:gamma_close2}).
Similarly, (\ref{eq:lam_close}) follows from (\ref{eq:lam_close2}).

The equation $\psi(z)=e^z - z - 1 = 0$ has the solution $z = 0$.
In order to show that it has a root $z \neq 0$, write
$z = x + iy$ to get
$$
e^x \cos y = 1 + x \mbox{ \rm and } e^x \sin y = y \, .
$$
Solve for $x$ to get
$x = y \cos y / \sin y - 1$. Inserting this value of $x$ into the second
equation we get
\begin{equation} \label{eq:y}
  \frac{y}{\sin y} = \exp \left(\frac{y \cos y}{\sin y} - 1 \right).
\end{equation}
We will find a solution of the form  $y = 2 \pi m + a$, where $\pi/4 < a < \pi/2$.
Note that if $y = 2 \pi m + \pi/4$, then the left hand side of (\ref{eq:y}) is
$\sqrt{2} y$, while the right hand side
is~$\exp(y-1)$, which is
strictly larger than $\sqrt{2}y$ for all $m\ge 1$.
If, on the other hand, $y = 2 \pi m + \pi/2$, then the left hand side
is~$y$ while the right hand side is~$\exp(-1)$, which is strictly
smaller than~$y$.
We  conclude that for all integers~$m\ge 1$, there exists
at least one solution   $y = 2 \pi m + a$, where $\pi/4 < a < \pi/2$.
\end{proof}

\subsection{The test function}
In this subsection we fix an eigenvalue $\lam$ of $M$ such that
$|\lam| \geq 1 - O(1/n)$, and let $f:[n]\to\C$ be a corresponding
eigenfunction.  We will denote the states of the $n$ cards at time $t$
by $\sigma_t(0),\ldots,\sigma_t(n-1)$, and assume that at time~$0$ we
start with the identity permutation, so $\sigma_0(i) = i$ for all $i$.
We emphasize that $\sigma_t$ is obtained from $\sigma_{t-1}$ by first
transposing the card at state $0$ with a uniform random card, and then
moving all cards one state up (modulo n). Thus for each $i$, the
sequence $\{\sigma_t(i)\}_{t \ge 0}$ is a Markov chain with transition
matrix $M$. To relate this to the description of the cyclic-to-random
shuffle $\{\sigma_t^*\}$ in the introduction, observe that
$\sigma_t^*$ is obtained from $\sigma_t$ by a rotation of size $t \!
\mod n$.

We will focus on the following test function $F:S_n \to \C$ :
% for the cyclic-to-random shuffle:
\begin{equation} \label{eq:defF}
F(\sigma) = \frac{1}{n} \sum_{i=0}^{n-1} f(\sigma(i)) \overline{f(i)}.
\end{equation}
Since  $f$ satisfies $\sum_{i=0}^{n-1} f(i) = 0$, under
the uniform distribution $\Unif$ on $S_n$ we have
\begin{equation} \label{eq:0stmoment}
\E_\Unif[F(\sigma)] = 0.
\end{equation}
It is also easy to see that for the  cyclic-to-random shuffle,
\begin{equation} \label{eq:1stmoment}
\E[F(\sigma_t)] = \frac{1}{n} \sum_{i=0}^{n-1}
\E[f(\sigma_t(i))]
\overline{f(i)}
= \lam^t F(\sigma_{0}) = \lam^t \|f\|_2^2,
\end{equation}
where $\|\cdot\|_2$ denotes the $\ell_2$-norm w.r.t.\ the uniform
distribution on $[n]$, i.e., $\|f\|_2^2 =
{\frac{1}{n}}\sum_{i=0}^{n-1}
 |f(i)|^2$.

We now calculate the second moment of $F(\sigma)$ under the stationary
distribution.
\begin{Lemma} \label{lem:2nd_stat}
\[
\E_{\Unif}\Bigl(|F(\sigma)|^2\Bigr) = \frac{\|f\|_2^4}{n-1} \,.
\]
\end{Lemma}

\begin{proof}
We have
\begin{equation} \label{eq:sec_moment_stat}
\E_{\Unif}\Bigl(|F(\sigma)|^2\Bigr) =
\frac{1}{n^2} \sum_{i\ne j} \E_{\Unif}
\Bigl(f(\sigma(i)) \overline{f(\sigma(j))}\Bigr)
f(j) \overline{f(i)} +
\frac{1}{n^2} \sum_{i} \E_{\Unif}\Bigl(|f(\sigma(i))|^2\Bigr) |f(i)|^2.
\end{equation}
The second term in~(\ref{eq:sec_moment_stat}) can be evaluated as
\begin{equation} \label{eq:sec_moment_diag}
\frac{1}{n^2} \sum_{i} \E_{\Unif}\Bigl(|f(\sigma(i))|^2\Bigr) |f(i)|^2
= \frac{\|f\|_2^2}{n^2} \sum_{i} |f(i)|^2
= \frac{\|f\|_2^4}{n}.
\end{equation}
Now let $i \neq j$ and let $\eta$ be an independent copy of $\sigma$. Then
\[
\E_{\Unif}[f(\sigma(i)) \overline{f(\sigma(j))}] =
\frac{n}{n-1}
\left( \E_{\Unif}\Bigl(f(\sigma(i)) \overline{f(\eta(j))}\Bigr) -
       \frac{1}{n}\E_{\Unif}\Bigl(|f(\sigma(i))|^2\Bigr) \right) =
-\frac{ \E_{\Unif}\Bigl(|f(\sigma(i))|^2\Bigr)}{n-1} = -\frac{\|f\|_2^2}{n-1}.
\]
Similarly,
\[
\sum_{i \neq j} f(j) \overline{f(i)} =
\sum_{i} \sum_{j} f(j) \overline{f(i)} - \sum_{i} |f(i)|^2
= -n \|f\|_2^2.
\]
Therefore the first term in~(\ref{eq:sec_moment_stat}) can be evaluated as
\begin{equation} \label{eq:sec_moment_non_diag}
\frac{1}{n^2} \sum_{i \neq j} \E_{\Unif}[f(\sigma(i)) \overline{f(\sigma(j))}]
f(j) \overline{f(i)} =
-\frac{\|f\|_2^2}{n^2 (n-1)} \sum_{i \neq j} f(j)
\overline{f(i)} =
\frac{\|f\|_2^4}{n(n-1)}.
\end{equation}
Combining~(\ref{eq:sec_moment_diag}) and~(\ref{eq:sec_moment_non_diag}),
we obtain the result via~(\ref{eq:sec_moment_stat}).
\end{proof}

For later use, we record here a simple variational bound on~$f$:
\begin{Lemma} \label{lem:2nd_inf}
There exists a universal constant~$c$ such that
$\|f\|_{\infty} \leq c \|f\|_2$ for all~$n$.
\end{Lemma}

\begin{proof}
Recall that $|\gamma - 1| \leq c/n$ for some $c > 1$ and
that $|\gamma| \leq 1$.
It follows from~(\ref{eq:x_i}) that, for all $k \neq 0$,
\begin{equation} \label{eq:f_above}
|f(k)| \leq 1 + \frac{c}{n} \sum_{j=1}^{k-1} |\gamma^{j-1}| \leq 1 + c.
\end{equation}
On the other hand, for $1 \leq k \leq n/2c$ we have
\begin{equation} \label{eq:f_below}
|f(k)| \geq 1 - \frac{c}{n} \sum_{j=1}^{k-1} |\gamma^{j-1}| \geq
1 - \frac{c k}{n} \geq \frac{1}{2}.
\end{equation}
By (\ref{eq:f_above}), we have $\|f\|_\infty \leq 1 + c$.
{}From (\ref{eq:f_below}), it follows that
\[
\|f\|_{2} \geq
\left( \frac{1}{2c} \left( \frac{1}{2} \right)^2 \right)^{1/2}
\geq \frac{1}{2} \frac{1}{(2c)^{1/2}}
\]
This completes the proof.
\end{proof}

\subsection{The second moment of  $F(\sigma_t)$.}
We begin with
an estimate of the contribution to the second moment from a specific
pair of cards.
Fix two distinct cards, $i$ and~$j$.
Denote by $A_{i}(s)=\{\sigma_s(i)=0\}$
the event that at step $s$ card $i$ is in state
$0$ (so it will be transposed with a uniform random card in the
next step). Let
$$N_{ij}(t)=\sum_{s=0}^{t-1} (\P[A_i(s)] + \P[A_j(s)])
$$
denote the expected number of times  $s<t$ where one of cards
$i,j$ was at state $0$.
Since at each step there is exactly one card at position $0$,
we have $\sum_{i=0}^{n-1} \sum_{s=0}^{t-1} \P[A_i(s)] = t$ and therefore
\begin{equation} \label{nij}
\sum_{i \neq j} N_{ij}(t) \leq  2n t \, .
\end{equation}

Next, we will couple~$\{\sigma_t\}$ with a
process $\{(\eta_t,\widetilde{\eta}_t)\}$, where $\eta$
and~$\widetilde{\eta}$ are two {\it independent\/} copies of the
cyclic-to-random shuffle starting from the identity permutation.  We
will observe the motions of cards~$i,j$ in $\eta,\widetilde{\eta}$
respectively; note that, unlike in~$\sigma$, these two motions are
independent.  We use the coupling to bound the dependence between the
cards in~$\sigma$.

\begin{Lemma} \label{lem:newcoupling}
For any two cards $i\ne j$ and all $t$ we have
\begin{equation}
\label{eq:coupling}
\left|\E\left[f(\sigma_t(i)) \overline{f(\sigma_t(j))}\right] -
\E\left[f(\eta_t(i))\overline{f(\widetilde{\eta}_t(j))}\right] \right| \leq
\frac{4 t + 4 n N_{ij}(t)}{n^2}\|f\|^2_{\infty}.
\end{equation}
\end{Lemma}

\begin{proof}
We define inductively a coupling of the process $\{\sigma_t\}$
and the pair process $\{(\eta_t,\widetilde{\eta}_t)\}$.
If $(\sigma_s(i),\sigma_s(j)) \neq (\eta_s(i),\widetilde{\eta}_s(j))$ then the
updates for the $\sigma$ and $(\eta,\widetilde{\eta})$ are performed
independently.
Otherwise, we have
\begin{equation} \label{sweet}
(\sigma_s(i),\sigma_s(j)) = (\eta_s(i),\widetilde{\eta}_s(j)) \,,
\end{equation}
and there are three cases to consider in the definition of the coupling at
step $s+1$:
\begin{description}
\item{{\bf Case 1.}}
Card $i$ is at position $0$ at time $s$.
\item{{\bf Case 2.}}
Card $j$ is at position $0$ at time $s$.
\item{{\bf Case 3.}}
Both of the cards $i,j$ are not in position $0$ at time $s$.
\end{description}
%The first two cases are treated similarly.
In Case 1, we take
%,  for all $\ell \neq  \sigma_s(j)+1$:
\[
(\sigma_{s+1}(i),\sigma_{s+1}(j))  = \left\{ \begin{array}{lll}
 (\ell,\sigma_s(j) + 1) & \mod n             & \mbox{w.p. } \frac{1}{n}
 \mbox{ \rm  for all } \ell \neq  \sigma_s(j)+1 ,  \\[1ex]
 (\sigma_s(j) + 1,1)  & \mod n             & \mbox{w.p. } \frac{1}{n}, \\
 \end{array} \right.
\]
and % for all $\ell \neq  \widetilde{\eta}_s(j)+1$:
\[
(\eta_{s+1}(i),\widetilde{\eta}_{s+1}(j)) = \left\{ \begin{array}{lll}
 (\ell,\widetilde{\eta}_s(j) + 1) & \mod n &
 \mbox{w.p. } \frac{n-1}{n^2}
 \mbox{ \rm for all } \ell \neq  \widetilde{\eta}_s(j)+1 \\[1ex]
 (\widetilde{\eta}_s(j) + 1,1)                   & \mod n &
 \mbox{w.p. } \frac{1}{n^2}, \\
 (\widetilde{\eta}_s(j) + 1,\widetilde{\eta}_s(j)+1)   & \mod n &
 \mbox{w.p. } \frac{n-1}{n^2}, \\[1ex]
 (\ell,1)                     & \mod n &
 \mbox{w.p. } \frac{1}{n^2}, \\[1ex]
 \end{array} \right.
\]
Thus, given that the processes satisfy (\ref{sweet})
 at time $s$ and that at that time
card $i$ is at location $0$, we may couple the processes
to satisfy (\ref{sweet}) at time
$s+1$ with conditional probability at least $\frac{(n-1)^2}{n^2}>
1-\frac{2}{n}$.
Similarly, in Case 2, if the coupling satisfies (\ref{sweet})
 at time $s$ then (\ref{sweet})  can be satisfied
at time  $s+1$ with conditional probability at least $1-\frac{2}{n}$.

In Case 3, the transition probabilities for the process $\sigma$
are given by
\[
(\sigma_{s+1}(i),\sigma_{s+1}(j))  = \left\{ \begin{array}{lll}
 (\sigma_s(i) + 1,\sigma_s(j) + 1) & \mod n & \mbox{w.p. } 1 - \frac{2}{n}, \\
 (\sigma_s(i) + 1,1)  & \mod n             & \mbox{w.p. } \frac{1}{n}, \\
 (1,\sigma_s(j) + 1)  & \mod n             & \mbox{w.p. } \frac{1}{n}. \\
 \end{array} \right.
\]
and the transition probabilities for the process $(\eta,\widetilde{\eta})$
are given by
\[
(\eta_{s+1}(i),\widetilde{\eta}_{s+1}(j)) = \left\{ \begin{array}{lll}
 (\eta_s(i) + 1,\widetilde{\eta}_s(j) + 1) & \mod n &
 \mbox{w.p. } 1 - \frac{2}{n} + \frac{1}{n^2}, \\
 (\eta_t(i) + 1,1)                   & \mod n &
 \mbox{w.p. } \frac{1}{n} - \frac{1}{n^2}, \\
 (1,\widetilde{\eta_s}(j) + 1)                     & \mod n &
 \mbox{w.p. } \frac{1}{n} - \frac{1}{n^2}, \\
 (1,1)                                  & \mod n &
 \mbox{w.p. } \frac{1}{n^2}.           \\
 \end{array} \right.
\]
It therefore follows that in Case 3, if the processes satisfy
(\ref{sweet}) at time $s$,
they may be coupled to satisfy it at time $s+1$ with conditional
probability  $1-\frac{4}{n^2}$.

It now follows that the probability that the processes ``unglue''
by time $t$ (i.e., (\ref{sweet}) fails for some $s \le t$)
is at most
\begin{equation} \label{eq:non_coupling}
\frac{2}{n} N_{ij}(t) +\frac{2t}{n^2} \,.
\end{equation}

We now estimate the difference of expected values in (\ref{eq:coupling}).
On the event where the processes satisfy (\ref{sweet}) at time $t$
 we get a $0$ contribution. On the
complementary event we get a contribution bounded by $2 \|f\|_{\infty}^2$.
We thus obtain (\ref{eq:coupling}) by (\ref{eq:non_coupling}).
\end{proof}

Since the processes $\eta$ and $\widetilde{\eta}$ are independent, it
follows from (\ref{eq:1stmoment}) that
\[
\E\left[f(\eta_t(i))\overline{f(\widetilde{\eta}_t(j))}\right] =
\E[f(\eta_t(i))]\E[\overline{f(\widetilde{\eta}_t(j))}] =
\lam^t f(i) \overline{\lam^t f(j)} = |\lam|^{2t} f(i) \overline{f(j)}.
\]
Therefore, by Lemma \ref{lem:newcoupling} we obtain
\begin{Corollary} \label{cor:newcoupling}
For any two cards $i\ne j$ we have
\begin{equation}
\label{eq:coupling_cor}
\Bigl|\E\Bigl(f(\sigma_t(i)) \overline{f(\sigma_t(j)}\Bigr)\Bigr|
\le
\left(|\lambda|^{2t}+ \frac{4 t + 4 n N_{ij}(t)}{n^2}\right)
\|f\|^2_{\infty}.
\end{equation}
\end{Corollary}
We can now bound the second moment of $F$.
\begin{Lemma} \label{lem:2nd_hard}
$ \displaystyle
\E\left[|F(\sigma_t)|^2\right] \leq
\left( |\lam|^{2t} + \frac{12 t + n}{n^2} \right) \|f\|_{\infty}^4.
$
\end{Lemma}
\begin{proof}
We have
\begin{equation} \label{eq:sec_moment_hard_dyn}
\E\left[|F(\sigma_t)|^2\right] =
\frac{1}{n^2} \sum_{i \neq j}
\E\left[f(\sigma_t(i)) \overline{f(\sigma_t(j))}\right] f(j) \overline{f(i)} +
\frac{1}{n^2} \sum_{i} \E\left[|f(\sigma_t(i))|^2\right] |f(i)|^2.
\end{equation}
Note that
\begin{equation}   \label{eq:sec_moment_hard_diag}
\frac{1}{n^2} \sum_{i} \E\left[|f(\sigma_t(i))|^2\right] |f(i)|^2
\leq  \frac{\|f\|_\infty^4}{n}.
\end{equation}
By Corollary \ref{cor:newcoupling}, for any $i\ne j$,
\begin{equation}  \label{eq:sec_moment_hard_non_diag}
\left|\E\left[f(\sigma_t(i)) \overline{f(\sigma_t(j))}\right] f(j)
\overline{f(i)}\right|
%&\leq& |\lam|^{2t} |f(i)|^2 |f(j)|^2 +
%\frac{4 t + 4 n N_{ij}(t)}
%{n^2} \|f\|_{\infty}^4 \\
\leq \left(|\lam|^{2t} +
\frac{4 t + 4 n N_{ij}(t)}{n^2} \right) \|f\|_{\infty}^4.
\end{equation}
Inserting~(\ref{eq:sec_moment_hard_diag})
and~(\ref{eq:sec_moment_hard_non_diag})
into~(\ref{eq:sec_moment_hard_dyn}) we  obtain
\begin{eqnarray*}
\E\left[|F(\sigma_t)|^2\right] \leq
\frac{\|f\|_{\infty}^4}{n^2}
\left(n + n^2 |\lam|^{2t} +
4t + \frac{4n}{n^2}  \sum_{i \neq j}
N_{ij}(t) \right)
%\\ &\leq&
%\frac{\|f\|_{\infty}^4}{n^2}
% \left(n + n^2 |\lam|^{2t} + 4t + \frac{4n}{n^2} 2nt \right)
\le \frac{\|f\|_{\infty}^4}{n^2} \left(n + n^2 |\lam|^{2t} + 12t \right),
\end{eqnarray*}
using (\ref{nij}). This completes the proof.
\end{proof}

\subsection{The mixing time}
Given the bound on the second moment of our test function
from the previous section, and the bound on the eigenvalue
from section~\ref{sec:roots}, it is straightforward to derive a lower
bound on the mixing time.

\par\noindent
{\bf Proof of Theorem \ref{th:lower}}
Recall from Lemma \ref{lem:delta_geq_1} that the equation $e^z-z-1=0$
has nonzero roots and let $\zeta$ be such a root.
By Lemma \ref{lem:delta_geq_1} it follows that for large $n$ there
exists a solution $\gamma$ of the equation
$(n-1) \gamma^n - n \gamma^{n-1} + 1 = 0$ satisfying
(\ref{eq:gamma_close}) and (\ref{eq:lam_close}). Fix $\gamma$ and let
$f$ be the corresponding eigen-function of $M$. Let $F$ be the test
function (\ref{eq:defF}). Write $\rho = n|1-\lam|$ and note that
$\rho = |\zeta+1|+O(\frac{1}{n})$ by (\ref{eq:lam_close}).

We use the test function $F$. Let $\mu_t$ be the distribution of $\sigma_t$ in the
cyclic-to-random shuffle where $\sigma_0$ is the identity permutation
and recall that $\Unif$ denotes the (uniform) stationary measure on $S_n$.
Let $g^2$ be the density of $\mu_t$
with respect to $\nu = (\mu_t + \Unif)/2$. Let $h^2$
be the density of $\Unif$
with respect to $\nu$.

By (\ref{eq:0stmoment}) and (\ref{eq:1stmoment}) we have that
\[
|\lam|^t \|f\|_2^2 = \left|\E_{\mu_t}[F] - \E_{\Unif}[F]\right| =
\left|\int F g^2 \, d\nu - \int F h^2 \, d\nu\right|.
\]
On the other hand, by Cauchy-Schwartz,
\[
\left|\int F g^2 \, d\nu - \int F h^2 \, d\nu\right|^2 =
\left|\int F (g+h)(g-h) \, d\nu\right|^2 \leq
\int |F|^2 (g + h)^2 \, d\nu
\cdot \int (g - h)^2 \, d\nu \,.
\]
By Lemma \ref{lem:2nd_stat} and Lemma \ref{lem:2nd_hard},
\begin{eqnarray*}
\int |F|^2 (g+h)^2 \, d\nu &\leq& 2 \int |F|^2 g^2
 \, d\nu + 2 \int |F|^2 h^2 \, d\nu
= 2\E_{\mu_t}\Bigl(|F|^2\Bigr) + 2\E_{\Unif}\Bigl(|F|^2\Bigr) \\
&\leq&
\frac{2 \|f\|_2^4}{n-1} +
2 \left( |\lam|^{2t} + \frac{12 t + n}{n^2} \right) \|f\|_{\infty}^4 \leq
2 \left( |\lam|^{2t} + \frac{12 t + 3n}{n^2} \right) \|f\|_{\infty}^4.
\end{eqnarray*}
Moreover,
\[
\int (g-h)^2 \, d\nu \leq \int |g^2 - h^2| \, d\nu =
2\TVD{\mu_t-\Unif}.
\]

Recalling Lemma \ref{lem:2nd_inf}, we  conclude that
\[
\TVD{\mu_t-\Unif} \geq \frac{|\lam|^{2t} \|f\|_2^4}
{4 \|f\|_{\infty}^4
\Bigl(|\lam|^{2t} + \frac{12 t + 3 n}{n^2}\bigr)} =
\Omega \left( \frac{|\lam|^{2t}}
{\Bigl(|\lam|^{2t} + \frac{12 t + 3n}{n^2}\Bigr)} \right) \,.
\]
It follows that $\TVD{\mu_t-\Unif} = \Omega(1)$ when
$\frac{12 t + 3n}{|\lam|^{2t} n^2} = O(1)$.
Note that if $t \geq n$, then
$\frac{12t + 3n}{|\lam|^{2t} n^2} = O(\frac{t}{n^2 |\lam|^{2t}})$
and that if
\[
t = \frac{1}{\rho} n \left(\log n - \log \log n -b \right) =
\frac{1+O(1/n)}{ |\zeta+1|} n \left(\log n - \log \log n -b \right) ,
\]
 then
\[
\frac{t}{n^2 |\lam|^{2t}} =
(1 + O(1/n))  \frac{n\log n }{\rho n  \log n e^b } =
  \frac{1+O(1/n) }{\rho  e^b } \,.
\]
The proof of the theorem follows.  \hfill\qedsymb

%%%%%%%%%%%%%%%%%%%%%%%%%%%%%%%%%%%%%%%%%%%%%%%%%%%%%%%%%%%%%%%%%%%%%%%%%%%
%%%%%%%%%%%%%%%%%%%%%%%%%%%%%%%%%%%%%%%%%%%%%%%%%%%%%%%%%%%%%%%%%%%%%%%%%%%
\section{An upper bound for general semi-random transpositions}
\label{sec:upper}
%%%%%%%%%%%%%%%%%%%%%%%%%%%%%%%%%%%%%%%%%%%%%%%%%%%%%%%%%%%%%%%%%%%%%%%%%%%
%%%%%%%%%%%%%%%%%%%%%%%%%%%%%%%%%%%%%%%%%%%%%%%%%%%%%%%%%%%%%%%%%%%%%%%%%%%

In this final section we prove Theorem~\ref{th:upper}.
\begin{proof}
By the triangle inequality it suffices to prove the theorem assuming
that the $L_t$ are deterministic. We thus restrict to that case.

We define a {\it strong uniform time\/} for the shuffle, i.e., a stopping
time~$T$ with the property that, given $T=t$, the random permutation
$\sigma^*_t$ has
the uniform distribution over~$S_n$.
It is well known (see, e.g.,~\cite{AlDi})
that, if $T$ is a strong uniform time, then the distribution
$\mu_t^*$ of $\sigma_t^*$ satisfies
$$
\TVD{\mu^*_t-\Unif} \le \P[T>t] \quad \forall t.
$$

Following Broder (as described in~\cite{DiaconisBook}) and Mironov~\cite{Mironov},
we define the stopping time in terms of a card marking process as follows.
Initially all cards are unmarked.
First, the card initially at $L_1$ is marked.
Later, at time~$t$, we mark the card at $L_t$ if it is unmarked and the card at
$R_t$ is already marked, and also if $R_t=L_t$  and this location
has an unmarked card.
Once a card is marked it remains so at all future times.  Set $T$ to be the
first time~$t$ at which all cards are marked.  Clearly $T$ is a stopping time.
The theorem follows immediately from the following two claims:
\par\noindent
{\bf Claim 1:} $T$ is a strong uniform time.
\par\noindent
{\bf Claim 2:} There exists  $C_0<\infty$ such that for any  $C_1 >
C_0$ we have
$$
\P\Bigl(T> C_1 n\log n \Bigr)\le n^{-\beta},
$$
for some $\beta=\beta(C_1)>0$.
Specifically, this holds for $C_0= 32\theta^{-3}+\theta^{-1}$,
where   $\theta =  e^{-2}(1-e^{-1})/2$.

\par\noindent
\noindent{\bf Proof of Claim~1:}
By induction, it is easy to check the following.
At any time $t$, given that $k$ cards have been marked,
conditional on the set of marked cards and their locations,
the mapping between these two sets (assigning to every marked card
its location) is uniformly distributed among the $k!$ possibilities.
See \cite{Mironov} or  \cite{DiaconisBook} for details.

\par\noindent
{\bf Proof of Claim 2:}
Divide time into successive {\it epochs\/} of length~$2n$,
starting after the card at $L_1$ is marked.
Denote by $u_k$ the fraction of unmarked cards before epoch $k$, so
$u_1=1-1/n$. Let $m_k=1-u_k$.  Let ${\cal H}_k$ denote the
history of the process prior to epoch~$k$,
and note that $u_k$ is a function of  ${\cal H}_k$.

\noindent{\bf Claim 3}: $\E(u_{k+1} | {\cal H}_k) \le u_k[1-2\theta m_k]$
for all $k$, where $\theta =  e^{-2}(1-e^{-1})/2$.
\begin{proof}
Consider a card $x$, unmarked before epoch $k$.
Of the $2n$ prescribed locations $\{L_t\}$ in the epoch, at most $n$
are  their last occurrence in the epoch. Thus for $1 \le j \le n$ we can find
$t(j)<s(j)$ in the epoch such that $L_{t(j)}=L_{s(j)}$.
For each $j \le n$, we have $R_{t(j)}=x$ with probability $1/n$.
Therefore, the event $A_x$ that there exists a $j \le n$ satisfying $R_{t(j)}=x$,
has probability
$$
\P(A_x |  {\cal H}_k ) \ge 1-(1-1/n)^n\ge 1-e^{-1}.
$$
On $A_x$, we fix $j$ to be minimal such that $R_{t(j)}=x$.  Given $A_x$ and
${\cal  H}_k$, with probability at least
$(1-1/n)^{2n-2}>e^{-2}$, we have $R_t \ne L_{t(j)}$ for all $t$ such that
$t(j) < t < s(j)$. In that case, $x$ is untouched by the random choices
between times $t(j)$ and $s(j)$, and then with probability at least $m_k$
the card at $R_{s(j)}$ is one of the $nm_k$ cards marked prior to epoch
$k$. Thus $x$ gets marked with probability at least $2\theta m_k$.
The assertion of Claim~3 follows.
\end{proof}

\par\noindent
{\bf Proof of Claim 2 continued:}
Using Claim~3, we first quantify the time to mark at least half the cards
(i.e., to achieve $m_k\ge 1/2$), and then the time to mark the remaining
cards (i.e., to achieve $u_k< 1/n$).
Denote by $D_k$ the number of cards that get marked
during epoch $k$ as a result of being transposed with a card that was marked
prior to epoch $k$. Clearly $m_{k+1}\ge m_k+D_k/n$.
The proof of Claim~3 implies that
\begin{itemize}
\item[(i)] if $m_k< 1/2$, then $\E(D_k | {\cal H}_k) \geq \theta n m_k$;
\item[(ii)] if $m_k \ge 1/2$,
then $\E(u_{k+1} | {\cal H}_k) \le (1-\theta)u_k$.
\end{itemize}

%Let $k_*$ denote the first $k$  for which $m_k\ge 1/2$.
%To bound the tail of $k_*$, we need a stochastic lower bound for $D_k$:
To bound the number of epochs where $m_k<1/2$,
we need a stochastic lower bound for  $D_k$:

\noindent{\bf Claim 4}: If $m_k < 1/2$, then
$$
\P\Bigl(D_k \ge \frac{\theta n m_k}{2}\,\Big|\,{\cal H}_k\Bigr) \ge \frac{\theta^2}{8}.
$$
\begin{proof}
Using the notation in the proof of Claim 3,
Denote by $\wdk$ the number of $j \le n$ such that $R_{s(j)}$
is one of the $nm_k$ cards marked prior to epoch $k$.
Clearly $D_k \le \wdk$.
The distribution of $\wdk$ is Binomial$(n,m_k)$, and this also
 holds given  ${\cal H}_k$. Therefore,
$$
\E(D_k^2 |  {\cal H}_k ) \le \E(\wdk^2 | {\cal H}_k )
\le (nm_k)^2+nm_k \le 2(nm_k)^2.
$$
In conjunction with (i) above, this yields
$$
\E(D_k^2 | {\cal H}_k) \le C_2 \E(D_k | {\cal H}_k)^2 \,,
$$
where $C_2=2\theta^{-2}$.
A standard second moment bound
(see, e.g., \cite[p.~8]{Kahane}) now yields Claim 4.
\end{proof}

\par\noindent
{\bf Proof of Claim 2 concluded:}
Call epoch~$k$ a ``growth epoch'' if
$m_{k+1} \ge (1+\theta/2) m_k $.
Call epoch~$k$ a ``good epoch'' if it is a growth epoch or
it satisfies $m_k \ge 1/2$.
 Claim~4 implies that the conditional probability that epoch $k$ is a good
 epoch, given ${\cal H}_k$, is at least~$\theta^2/8$. Thus the number of
 good epochs among the first $k_3=C_3 \log n$ epochs stochastically
 dominates a Binomial$(k_3, \theta^2/8)$ random variable.
Fix $C_3> 32\theta^{-3}$, and
denote by $\Omega_3$ the event that there are at least $(4\log n) /\theta$ good
epochs among the first $k_3$ epochs.
Recall that the probability that a binomial random variable differs from its
mean by a constant multiple of the mean decays exponentially in the
number of trials $k_3=C_3 \log n$.
We infer that  $\P(\Omega_3^c) <n^{-\beta}/2$
for some $\beta>0$.
Since $(1+\theta/2)^{4/\theta}>e$ and $m_1=1/n$,
the number of growth epochs must be smaller than $(4\log n)/\theta$.
Thus on $\Omega_3$ we have $m_{k_3} \ge 1/2$.

\medskip

Turning now to the second portion, once $m_k\ge 1/2$ we have from~(ii)
above that $\E(u_{k+1} | u_k) \le (1-\theta) u_k$.
Therefore for all $k>0$, we have
$$
\E(u_{k_3+k}\, |\, \Omega_3, u_{k_3}) \le (1-\theta)^k u_{k_3} \le
 \frac{e^{-\theta k}}{2}\,.
$$
Thus if $k=(1+\beta)\theta^{-1}\log n$, then
$$ \P(u_{k_3+k} \ge 1/n\,|\, \Omega_3) \le \E(n u_{k_3+k}\, |\, \Omega_3)
   \le n\cdot \frac{n^{-1-\beta}}{2}=\frac{n^{-\beta}}{2} \,.
$$
In conjunction with the bound for $\P(\Omega_3^c)$,
this implies that $ \P(u_{k_3+k} \ge 1/n) \le n^{-\beta}$
for this value of $k$.
In other words, if $C_1 > C_0= 32\theta^{-3}+\theta^{-1}$
and $k_1=C_1 \log n$, then there exists $\beta=\beta(C_1)>0$ such that
$\P(u_{k_1} \ge 1/n) \le n^{-\beta}$.
This completes the proofs of Claim~2 and of the theorem.
\end{proof}

\section{Concluding remarks and further problems}
\begin{enumerate}
\item We have shown that the cyclic-to-random transposition shuffle
on $n$ cards has mixing time of order $\Theta(n \log n)$.
However, the constant in our general upper bound, and that in the specific
cyclic-to-random upper bound of Mironov~\cite{Mironov}, are significantly
larger than the constant in our lower bound.
We believe that the lower bound is closer to the truth,
and moreover, that this shuffle exhibits the ``cutoff phenomenon'',
i.e., there is a constant $C_*$ such that for $t<(1-\epsilon)C_* n \log
n$ the distribution after $t$ steps, $\mu^*_t$, satisfies
$\TVD{\mu^*_t-\Unif} =1-o(1)$ as $ n \to \infty$, while for
$t>(1+\epsilon)C_* n \log n$ we have $\TVD{\mu^*_t-\Unif} =o(1)$ as $
n \to \infty$. Proving this, and determining $C_*$, remain a challenge.

\item
Does the cyclic-to-random shuffle capture the key features
of the RC4 cryptographic system, as suggested by
Mironov~\cite{Mironov}? \newline
If the answer is positive, then we expect that the test function
$F$ defined in  (\ref{eq:defF}) may play a role in future analysis of RC4.

\item For which sequence $\{L_t\}$ does the resulting
semi-random transposition shuffle on $n$ cards have the largest mixing
time?
\newline
We suspect that the slowest shuffle in this class is the "star transpositions"
shuffle, for which $L_t=0$ for all $t$, and the mixing time is $(1+o(1)) n \log n$
by \cite{diaconis2}.

\item
Is there a universal constant $c>0$ such that, for any
semi-random transposition shuffle on $n$ cards,
the mixing time is at least $c n \log n$ ? \newline
For this lower bound question there is no obvious reduction to the
case where the sequence $\{L_t\}$ is deterministic, so conceivably
the question could have different answers for
deterministic  $\{L_t\}$ and random  $\{L_t\}$.
Two specific cases of interest are:
\begin{itemize}
\item For each $k \ge 0$, let $\{L_{kn+r}\}_{r=1}^n$ be a uniform random
  permutation of $\{0,\ldots,n-1\}$, where these permutations are
independent.
 \item Let $\{L_t\}$ be a Markov chain with memory 2, where
$L_1=0,L_2=1$ and for each $t \ge 3$ we have $L_{t+1}=2L_t-L_{t-1} \!
\mod n$ with probability $1-1/n$ and $L_{t+1}=L_{t-1} $
with probability $1/n$. This choice of  $\{L_t\}$ was suggested to us
by Igor Pak (personal communication), motivated by \cite{CLP}.
\end{itemize}
Each of these examples has a ``quenched'' version, where the sequence
 $\{L_t\}$ is picked in advance and then used as a deterministic
 sequence, and an ``annealed'' version, where
the $\{L_t\}$ are random variables with the specified distribution.

\end{enumerate}

\noindent{\bf Acknowledgments:} We are grateful to Serban  Nacu for his help
with some of the complex analysis used in this paper.

%%%%%%%%%%%%%%%%%%%%%%%%%%%%%%%%%%%%%%%%%%%%%%%%%%%%%%%%%%%%%%%%%%%
%%%%%%%%%%%%%%%%%%%%%%%%%%%%%%%%%%%%%%%%%%%%%%%%%%%%%%%%%%%%%%%%%%%

\end{document}